\documentclass[final]{article}%
\usepackage{isipta2009}
\usepackage{amsmath}
\usepackage{amsfonts}
\usepackage{amssymb}
\usepackage{cite}
\usepackage{graphicx}
\usepackage{xcolor}
\usepackage{pgf}%
\newtheorem{theorem}{Theorem}

\newtheorem{definition}{Definition}
\newtheorem{example}{Example}

\newenvironment{proof}[1][Proof]{\noindent\textbf{#1.} }{\ \rule{0.5em}{0.5em}}
\begin{document}

\title{Dutch Books and Combinatorial Games}
\author{Peter Harremo\"{e}s\\Centrum Wiskunde \& Informatica (CWI), The Netherlands\\P.Harremoes@cwi.nl }
\maketitle

\begin{abstract}
The theory of combinatorial games (like board games) and the theory of social
games (where one looks for Nash equilibria) are normally considered two
separate theories. Here we shall see what comes out of combining the ideas. J.
Conway observed that there is a one-to-one correspondence between the real
numbers and a special type of combinatorial games. Therefore the payoffs of a
social games are combinatorial games. Probability theory should be considered
a safety net that prevents inconsistent decisions via the Dutch Book Argument.
This result can be extended to situations where the payoff function yields a
more general game than a real number. The main difference between
number-valued payoff and game-valued payoff is that the existence of a
probability distribution that gives non-negative mean payoff does not ensure
that the game will not be lost.

\end{abstract}

\begin{keywords}
Combinatorial game, Dutch Book Theorem, exchangable sequences, game theory, surreal number.
\end{keywords}

\section{Introduction}

The word game in mathematics has two different meanings. The first type of
games are the \emph{social games} where a number of agents at the same time
have to make a choice and where the payoff to each agent is a function of all
agents' choices. Each agent has his own payoff function. The question is how
the agents should choose in order to maximize their own payoff. In general the
players may benefit by making coalitions against each other. This kind of game
theory has found important applications in social sciences and economy. A
special class of these social games are the two-person zero-sum games where
collaboration between the agents makes no sense.

The second type of games are the \emph{combinatorial games}. These are
mathematical models of board games. These games are the ones that people find
interesting and amusing. Games that people play for amusement often involve an
element of chance, generated by, for instance, dice, but the combinatorial
games are by definition the ones that do not contain this element. Therefore
they are sometimes called \emph{games of no chance} \cite{Nowakowski1996}.
Examples from this category are chess, nim, nine-mens-morris, and go.
Combinatorial game theory has been particularly successful in the analysis of
impartial games like nim \cite{Berlekamp1982} and has lead to a better
understanding of endgames in go \cite{Berlekamp1991, Berlekamp1996,
Nowakowski1996}.

The Dutch Book Theorem is important in our understanding of imprecise
probabilities. The Dutch Book Theorem was first formulated and proved by F. P.
Ramsay in 1926 (reprinted in \cite{Ramsey1931}) and later independently by B.
de Finetti \cite{Finetti1937}, who used it as an argument for a subjective
interpretation of probabilities. Since the original formulation of the Dutch
Book Theorem most of the research has been in the direction of more subjective
versions. As it is normally formulated, the theorem relies on the concept of a
\emph{real-valued payoff} function. One may think of an outcome of the payoff
function as money but the uniform mean of having \pounds \ 1.000.000\ and
having \pounds \ 0\ is having \pounds \ 500.000. Most people have a very clear
preference for having \pounds \ 500.000\ rather than an unknown amount of
money with mean \pounds \ 500.000. Instead one may think of the payoff as a
more subjective notion of \emph{value, }but this is also a highly debatable
concept and one may actually consider money as our best attempt to quantify
value. Savage showed that the concept of value and payoff function can be
replaced by the concept of preference, so that a coherent set of preferences
corresponds to the existence of a payoff function and a probability measure.
This line of research has been followed up by many other researchers
\cite{Savage1954, Chu2008}. All those studies involve some subjective notion
of value or preference.

In order to better understand the Dutch Book Theorem it is desirable to see
how the theory would look in an environment where a subjective notion of value
plays no role. In this study we replace the normal payoff functions by
game-valued functions. There are several reasons why this is of interest:

\begin{itemize}
\item A real-valued payoff function is a special case of a game-valued payoff function.

\item The theory of probability has its origin in the study of games involving chance.

\item Social game theory and combinatorial game theory may mutually benefit
from a closer interaction.

\item One can often get insight into a special case by the study of its generalizations.
\end{itemize}

With a game-valued payoff function the players in a social game have to play a
certain combinatorial game that depends on their decisions and/or on some
random event. This setup may seem quite contrived, but many board games that
involve chance are of this form.

\begin{example}
In chess it is normally considered a slight advantage to play white. Therefore
one normally randomly selects who should play white and who should play black.
\end{example}

\begin{example}
M. Ettinger has developed an interesting version of combinatorial game theory
where after each move a coin is flipped to determine who is going to play next
\cite{Ettinger1996}.
\end{example}

Actually any board game involving chance may be considered as an example. It
will be the subject of a future paper how to take advantage of a combined
probabilistic and combinatorial game approach for some specific board games.
In this short note we shall focus entirely on how we should formulate or
reformulate the Dutch Book Theorem when the payoffs are combinatorial games.

Social games and combinatorial games are built on quite different ideas and
many scientists only know one of the types of game theory. There have only
been few attempts to combine the two types of game theory \cite{Zhao2006,
Ettinger1996}. In this exposition we will assume that the reader has basic
knowledge about social games such as two-person zero-sum games. Nevertheless
we have to repeat some of the elementary definitions from social game theory
in order to fix notation and, in particular, to avoid confusion with similar
but slightly different concepts from combinatorial game theory.

Our main result is that it is possible to formulate versions of the Dutch Book
Theorem for game-valued payoff functions, but there will be some important
modifications of the theorem. For instance our probability distributions will
not always be real-valued. In our approach the focus is on order structure
(induced by games) and its relation to decision theory. A somewhat orthogonal
approach was taken in \cite{Khrennikov1994} where the probabilities were
elements of a metric space with no order structure.

\section{Combinatorial games}

The theory of combinatorial games was developed by J. Conway as a tool to
analyze board games \cite{Conway1976, Berlekamp1982}. A short and more careful
exposition can be found in \cite{Schleicher2006}. In a board game the players
\emph{alternate} in making moves. Each move changes the configuration of the
pieces on the board to some other configuration but only certain changes are
allowed. It is convenient to call the two players \emph{Left} and
\emph{Right}. We shall often consider different board configurations as
different games. If $G$ denotes a game, i.e. a certain configuration then the
game is specified by the configurations $G^{L}$ that Left is allow to move to
and the configurations $G^{R}$ that Right is allowed to move to, and we write
$G=\left\{  G^{L}\mid G^{R}\right\}  .$ Note that we have not told who is
playing first, and therefore we have to describe it from both Left's and
Right's perspective. Now the point is that $G^{L}$ and $G^{R}$ are sets of
games, so a game is formally a specification of two sets of games. In a board
game it is nice to have many options to choose among and bad if there are only
few options. The worst case for Left is if there are \emph{no options left}
and in this case we say that Left has \emph{lost} the game. So Left has lost
the game if he is to move next and $G^{L}$ is empty. Similar Right loses the
game if it is Right to move and $G^{R}$ is empty. The rules of many board
games can be modelled in this way.

\input{day1.TpX}

\begin{example}
[Games illustrated in Figure \ref{day1}.]\label{Exday1}The game $\left\{
\varnothing\mid\varnothing\right\}  $ is a boring one. The one to move first
loses this game. This game is denoted $0.$

The game $\left\{  \varnothing\mid0\right\}  $ is lost by Left if Left has to
move first. If Right goes first Right has to choose $0.$ Now it is Left to
move but this is a losing position for the one who is going to move, so poor
Left loses. Thus Right always wins the game $\left\{  \varnothing
\mid0\right\}  .$ This game is denoted $-1.$

The game $\left\{  0\mid\varnothing\right\}  $ is lost by Right if Right has
to move first. If Left goes first Left has to choose $0.$ Now it is Right to
move but this is a losing position for the one who is going to move, so now
Left is happy again because he wins. Thus Left always wins the game $\left\{
0\mid\varnothing\right\}  .$ This game is denoted $1.$

Similarly we see that $\left\{  0\mid0\right\}  $ is won by the player that
moves first. This game is called \emph{star} and is denoted $\mathord\ast.$ In
Japanese go terminology such a position is called \emph{dame}.
\end{example}

Here we shall use the following recursive definition of a game.

\begin{definition}
A \emph{game} is a pair $\left\{  G^{L}\mid G^{R}\right\}  $ where $G^{L}$ and
$G^{R}$ are sets of already defined games.
\end{definition}

The \emph{status} of a game $G$ can be classified according to who wins if
both players play optimally. We define%
\[%
\begin{array}
[c]{cl}%
G=0, & \text{if second player wins;}\\
G<0, & \text{if Right wins whoever plays first;}\\
G>0, & \text{if Left wins whoever plays first;}\\
G\,\mathrel\Vert\,0, & \text{if first player wins.}%
\end{array}
\]

For a game $G$ we can reverse the role of Left and Right and call this the
\emph{negative of the game}. Formally we use the following recursive
definition.%

\input{mean0.TpX}

Left and Right can play two games in parallel. In every round each player
should make a move in one of the games of his own choice. Perhaps there are
urgent moves to be made in both games so the players have to prioritize in
which game it is most important to make the move. Several games played in
parallel is called the \emph{sum of the games}, and many positions in actual
board games can be understood as sums of sub-games. Combinatorial game theory
is essentially the theory of how to prioritize your moves in a board game that
has the structure of a sum of independent sub-games. Formally the sum of the
games $G$ and $H$ is defined recursively by%
\begin{multline*}
G+H=\\
\left\{  \left(  G^{L}+H\right)  \cup\left(  G+H^{L}\right)  \mid\left(
G^{R}+H\right)  \cup\left(  G+H^{R}\right)  \right\}  .
\end{multline*}
The sum of games is normally illustrated by the disjoint union of the game
trees of the individual games. The game $G-H$ is by definition the game
$G+\left(  -H\right)  .$

Now, we are able to define what it should mean that two games are equal. We
write $G=H$ if $G-H=0,$ i.e. second player wins $G-H.$ One can define $G>H$,
$G<H,$ and $G\mathrel\Vert H$ in the same way. We say that $G$ and $H$ are
\emph{confused} if $G\mathrel\Vert H.$ One can prove that $G=H$ if and only if
$G+K$ and $H+K$ have the same status for any game $K.$

With these operations the class of games has the structure as a partially
ordered Abelian group. Any Abelian group is a module over the ring of integers
with multiplication defined as follows. If $n$ is a natural number we define
$n\cdot G$ by%
\[
\overset{n\text{ times}}{\overbrace{G+G+\cdots+G}}~.
\]
If $n=0$ then $0\cdot G$ is by definition equal to $0.$ If $n$ is a negative
integer we define $n\cdot G$ to be equal to $\left(  -n\right)  \cdot\left(
-G\right)  .$

The equation $2\cdot G=0$ has $G=0$ as solution, but $G=\mathord\ast$ is also
a solution. Therefore there is in general no unique way of defining
multiplication of a game by $1/2,$ and the same holds for other non-integers.
From this point of view it is surprising that all dyadic fractions (rational
numbers of the form $n/2^{m}$) can be identified with games. One way of doing
it goes as follows.

\section{Numbers may be identified with games}

J. Conway discovered that all real numbers can be identified with games but
his construction will lead to a larger class of numbers called the
\emph{surreal numbers} (or Conway numbers). The surreal numbers were first
described in a mathematical novel by D. Knuth \cite{Knuth1974}, and later in
much detail by J. Conway \cite{Conway1976}. For newer and more complete
descriptions we refer to \cite{Gonshor1986, Alling1987}.

We have already defined the game $1$ so the integer $n$ is identified with the
game $n\cdot1.$ The game $\left\{  0\mid1\right\}  $ satisfies
\[
2\cdot\left\{  0\mid1\right\}  =1.
\]
Hence the $2^{-1}$ can be identified with the game $\left\{  0\mid1\right\}
.$ In general the game $\left\{  0\mid2^{-m}\right\}  $ satisfies%
\[
2\cdot\left\{  0\mid2^{-m}\right\}  =2^{-m}%
\]
so the fraction $2^{-\left(  m+1\right)  }$ can be identified with the game

\input{fractions.TpX}

Thus the fraction $n/2^{m}$ can be identified with the game $n\cdot2^{-m}.$ In
this way any dyadic fraction can be identified with a game.

A real number can be identified with a Dedekind section in the group of dyadic
fractions. In other words, a real number $r,$ can be identified with the
partition of the dyadic fractions into the sets
\begin{align*}
A  &  =\left\{  n\cdot1/2^{m}<r\mid m,n\in\mathbb{N}\right\}  ,\\
B  &  =\left\{  n\cdot1/2^{m}>r\mid m,n\in\mathbb{N}\right\}  .
\end{align*}
Now, $A$ and $B$ can be identified with sets of games and therefore $\left\{
A\mid B\right\}  $ is a game. When $r$ is a real number that is not a dyadic
fraction, it can be identified with the game $\left\{  A\mid B\right\}  .$ At
this step one has to check that the structure of the real numbers as an
ordered group is preserved under the embedding but this turn out to be the
case \cite{Conway1976}.

We have seen that real numbers may be identified with games, but combining the
definition of a game with the idea of a Dedekind section leads to the much
larger class of numbers called the \emph{surreal numbers}. Formally a surreal
number is a game of the form $\left\{  A\mid B\right\}  $ where $A$ and $B$
are sets of (already constructed) surreal numbers such that $a<b$ for $a\in A$
and $b\in B.$ That means that a surreal number can always be played as a
combinatorial game.

\begin{example}
The first transfinite ordinal number $\omega$ is identified with the game
$\left\{  \mathbb{N}\mid\varnothing\right\}  .$ The equation $\omega-\omega=0$
makes no sense in Cantor's arithmetic for transfinite ordinals or cardinals,
but if we identify $\omega$ with a game the equation makes sense, because we
have
\[
\omega-\omega=\left\{  1,2,3,\cdots\mid\varnothing\right\}  +\left\{
\varnothing\mid-1,-2,-3,\cdots\right\}  .
\]
This game is essentially like "my father has more money than your father" and
most children soon experience that one should not start in such a game. It is
clear that $\omega$ should not be interpreted as an amount but is better
understood as a huge set of options. Conway identified all Cantor's ordinal
numbers with surreal numbers, but Cantor and Conway use \emph{different}
additive structures so the identification is somewhat problematic. For
instance Conway's addition is commutative but Cantor's addition of ordinal
numbers is not. Here we shall use $\omega$ as a symbol for a game rather than
an ordinal in Cantor's sense.
\end{example}

Formally the surreal numbers are constructed by (transfinite) recursion. It
starts with the number $0=\left\{  \varnothing\mid\varnothing\right\}  .$ In
each recursion step one adds new surreal numbers to the ones already
constructed. Addition and multiplication extend to surreal numbers and with
these operations the surreal numbers are a maximal ordered field. Although the
definition of surreal multiplication is relevant for the next two sections we
cannot present the definition in this short note but have to refer to
\cite{Schleicher2006, Conway1976}. For most computations surreal numbers are
not different from real numbers but the topology is different.

A game $G$ is said to be \emph{infinitesimal} if $-2^{-m}\leq G\leq2^{-m}$ for
all natural numbers $m.$ The number $1/\omega$ is an example of an
infinitesimal number that is positive. Between any two different real numbers
there are more than continuously many surreal numbers, and the intersection of
the intervals $\left[  -2^{-m};2^{-m}\right]  $ contains infinitely many
\emph{infinitesimal numbers}. Formally there are so many surreal numbers that
they do not form a set but a class.

\section{Surreal probabilities and payoffs\label{SecSurreal}}

Here we will introduce a version of the \emph{Dutch Book Theorem} for surreal
payoff functions. Because of the somewhat different topology of the surreal
numbers, we have to be a little careful in the formulation and proof of the
Dutch Book Theorem. In particular some of the standard methods for proving
these results like the Hahn-Banach theorem and the separation theorem for
convex sets, do not hold in their normal formulation when we are using surreal
numbers. Those used to to non-standard analysis may note that what we are
doing is essentially to veify that our result may be formulated in first order language.

The setup is as follows. Alice wishes to make a bet on an outcome $a\in A$. A
bookmaker $b\in B$ offers the surreal payoff $g\left(  a,b\right)  $ (positive
or negative) if the outcome of a random event is $a\in A.$ Thus $\left(
a,b\right)  \rightarrow g\left(  a,b\right)  $ can be considered as a matrix
when $A$ and $B$ are finite sets. Alice should reject to play with a bookmaker
$b$ if Alice thinks that the payoff function $a\rightarrow g\left(
a,b\right)  $ is not favorable. For simplicity we shall assume that Alice
accepts the payoff functions offered by the bookmakers $b\in B.$ We recall
that a surreal number is a game so if the outcome is $a$ and the bookmaker is
$b$ then Alice has to play the game $g\left(  a,b\right)  $ against the
bookmaker with Alice playing Left and the bookmaker playing Right.

By a \emph{portfolio} we shall mean a probability vector $Q=\left(
q_{b}\right)  _{b\in B}$ on $B.$ In this section will allow the portfolio to
have surreal values. Such a portfolio is described by the payoff function%
\begin{equation}
a\rightarrow\sum_{b\in B}q_{b}\cdot g\left(  a,b\right)  , \label{DutchNeg}%
\end{equation}
A \emph{Dutch book} is a portfolio such that (\ref{DutchNeg}) is negative for
all $a\in A,$ i.e. the portfolio game will be lost by Alice for any value of
$a\in A$.

We assume that one of the bookmakers $b_{0}$ offers a payoff function
$g\left(  a,b_{0}\right)  =0$ for all $a\in A$ ($b_{0}$ acts like a bank with
interest rate 0). Let $Q$ be a portfolio and assume that there exists a Dutch
book $Q^{\prime}.$ If $Q$ has $B$ as support then $q_{\min}=\min_{b\in B}%
q_{b}>0$ and the payoff is%
\begin{multline*}
\sum_{b\in B}q_{b}\cdot g\left(  \cdot,b\right)  =\\
\sum_{b\in B}\left(  q_{b}-q_{\min}\cdot q_{b}^{\prime}\right)  \cdot g\left(
\cdot,b\right)  +q_{\min}q_{b}\sum_{b\in B}q_{b}^{\prime}\cdot g\left(
\cdot,b\right)  <\\
\sum_{b\in B}\left(  q_{b}-q_{\min}\cdot q_{b}^{\prime}\right)  \cdot g\left(
\cdot,b\right)  +\left(  q_{\min}\sum_{b\in B}q_{b}^{\prime}\right)  \cdot
g\left(  \cdot,b_{0}\right)  .
\end{multline*}
Hence Alice should reject to play with at least one of the bookmakers. If no
Dutch book exists the set of payoff functions is said to be \emph{coherent}.
The notion of convexity will be used, and in this section we allow surreal
coefficients in convex combinations.

\begin{theorem}
\label{DutchBookThm1}%
\index{Theorem for the Alternative for Matrices}%
Let $A$ and $B$ denote finite sets and let $\left(  a,b\right)  \rightarrow
g\left(  a,b\right)  $ denote a surreal valued payoff function. If the payoff
function is coherent then there exists non-negative surreal numbers $p_{a}$
such that $\sum p_{a}=1$ and
\begin{equation}
\sum_{a\in A}p_{a}\cdot g\left(  a,b\right)  \geq0 \label{Positiv}%
\end{equation}
for all $b\in B.$
\end{theorem}

\begin{proof}
Assume that $A$ has $d$ elements. Then each function $g\left(  \cdot,b\right)
$ may be identified with a $d$-dimensional surreal vector. Let $K$ be the
convex hull of $\left\{  g\left(  \cdot,b\right)  \mid b\in B\right\}  ,$ and
let $L$ denote the strictly negative surreal functions on $A.$ They are convex classes.

If $K$ and $L$ intersect then there exists non-negative surreal numbers
$q_{b}$ such that $\sum q_{b}=1$ and such that (\ref{DutchNeg}) defines a
strictly negative function.

Assume that $K$ and $L$ are disjoint. Then define $C=K-L$ as the class of
vectors $\bar{x}-\bar{y}$ where $\bar{x}$ in $K~$and $\bar{y}$ in $L.$ This is
convex and does not contain $\bar{0}.$ Now, $K$ is a polytope (convex hull of
finitely many extreme points) and $L$ is polyhedral (given by finitely many
inequalities), so $C$ is polyhedral. Hence, each of the faces of $C$ is given
by a linear inequality of the form $\sum_{a\in A}p_{a}\cdot g\left(  a\right)
\geq c$ for $g\in C$. The delta function $\delta_{\alpha}$ is non-negative so
if $g$ is in $C$ then $g-\ell\cdot\delta_{\alpha}$ is also in $C$ for $\ell$
positive. In particular
\begin{align*}
c  &  \leq\sum_{a\in A}p_{a}\cdot\left(  g-\ell\cdot\delta_{\alpha}\right)
\left(  a\right) \\
&  =\sum_{a\in A}p_{a}\cdot g\left(  a\right)  -\sum_{a\in A}p_{a}\ell
\delta_{\alpha}\left(  a\right) \\
&  =\sum_{a\in A}p_{a}\cdot g\left(  a\right)  -\ell\cdot p_{\alpha}%
\end{align*}
for all positive $\ell.$ Hence $p_{\alpha}\geq0$ for all $\alpha\in A.$
Further we know that $\bar{0}$ is not in $C$ so that $\sum_{a\in A}p_{a}%
\cdot0\geq c$ does not hold and therefore $c>0.$ In particular $p_{a}$ cannot
be 0 for all $a.$ The result follows by replacing $p_{a}$ by%
\[
\frac{p_{a}}{\sum_{a\in A}p_{a}}.
\]

\end{proof}

Note that our surreal valued version Dutch Book Theorem states there are
\emph{two exclusive }cases:

\begin{enumerate}
\item Dutch book.

\item Non-negative mean value.
\end{enumerate}

The theorem leads to surreal probabilities $p_{a}\geq0$. Due to the
normalization we do not have infinite probabilities, but there is no problem
in having infinitesimal probabilities. In general the probability distribution
will not be uniquely determined, but will merely be located in a non-empty
convex set (credal set). Therefore the Dutch Book Theorem suggests that
uncertainty about some unknown event should be represented by a \emph{convex
set of surreal probability distributions} rather than a single real valued
distribution. Real functions are special cases of surreal functions so even if
the payoff functions are real valued one can model our uncertainty by a convex
set of surreal probability distributions.

If either $g$ is acceptable or $-g$ is acceptable then it is called a
two-sided bet. In this case the convex set of probability distributions
reduces to a point. The term one-sided bet is taken from F. Hampel
\cite{Hampel1999}. In general people will find it difficult to decide that
either $g$ or $-g$ is acceptable and thus the two-sided bet is not realistic.
In De Finetti \cite{Finetti1937} only two-sided bets were considered. In our
formulation of the Dutch Book Theorem we just have a one-sided bet with a set
of acceptable payoff functions.

A special case that has been studied in great detail is when the functions
$g\left(  \cdot,b\right)  $ only assume two different values, i.e. $g\left(
\cdot,b\right)  $ has the form
\[
g\left(  a,b\right)  =\left\{
\begin{array}
[c]{cc}%
g_{1}\left(  b\right)  , & \text{for }a\in A_{b};\\
g_{2}\left(  b\right)  , & \text{for }a\notin A_{b}.
\end{array}
\right.
\]
Without loss of generality we may assume that $g_{1}\left(  b\right)
\geq0>g_{2}\left(  b\right)  .$ Then the $g$ is accepted when $P\left(
A_{b}\right)  g_{1}\left(  b\right)  +\left(  1-P\left(  A_{b}\right)
\right)  g_{2}\left(  b\right)  \geq0$ or equivalently
\begin{equation}
P\left(  A_{b}\right)  \geq\frac{-g_{2}\left(  b\right)  }{g_{1}\left(
b\right)  -g_{2}\left(  b\right)  }. \label{nedre}%
\end{equation}
We then define the \emph{lower provision function} \cite{Walley1991} by
\[
L\left(  A\right)  =\min P\left(  A\right)
\]
where the minimum is taken over all distributions $P$ that satisfies
(\ref{nedre}) for all $b\in B.$ One may form surreal lower provisions in the
same way as ordinary lower provisions are formed.

In this section we have seen that uncertainty may be identified with a convex
set of surreal-valued probability distribution, but often such convex sets
contain a lot of real-valued distributions. One may therefore ask whether the
surreal-valued distributions add anything to the theory. Are they of any use?
This we will try to answer in the next section.

\section{Two-person zero-sum games}

The theory of two persons zero sum games was founded by J. von Neumann%
\index{Neumann@von Neumann}
together with O. Morgenstern \cite{Neumann1947} and has been extended to
social games with more players. The readers who are interested in a deeper
understanding of the theory of social games should consult \cite{Straffin1993}
for an easy introduction or \cite{Gibbons1992} for a more detailed exposition.

A social game with $2$ players, that we will call Alice and Bob, is described
by $2$ sets of \emph{strategies }$A,B$ such that Alice can choose a strategy
from $A$ and Bob can choose a strategy from $B.$ If Alice choose $a$ and Bob
choose $b$ then the payoff for Alice will be $g\left(  a,b\right)  $ and the
payoff for Bob will be $-g\left(  a,b\right)  ,$ where $g$ is a function from
$A\times B$ to surreal numbers$.$ Alice and Bob will never collaborate in a
zero-sum game because what is good for one of the players is equally bad for
the other.

A pair of strategies $\left(  a,b\right)  $ is called a \emph{Nash
equilibrium}%
\index{Nash!equilibrium|see{Nash pair}}
if no player will benefit by changing his own strategy if the other player
leaves his strategy unchanged. If a game has a unique Nash pair%
\index{Nash!pair}
and both players are \emph{rational}, then both players should play according
to the Nash equilibrium%
\index{Nash pair}%
.

Assume that the players are allowed to use mixed strategies, i.e. choose
independent probability distributions over the strategies. The probabilities
are allowed to take surreal values. Let $P$ be the mixed strategy of Alice and
$Q$ be a mixed strategy of Bob. Then the \emph{mean payoff} for Alice is%
\[
g\left(  P,Q\right)  =\sum_{\left(  a,b\right)  }g\left(  a,b\right)  \cdot
p_{a}q_{b}.
\]
This number is considered as the payoff of the social game where mixed
strategies are allowed.

\begin{theorem}
\label{ZeroSumThm}Consider a game with surreal valued payoffs. If the players
are allowed to use mixed strategies, then the game has a Nash equilibrium%
\index{Nash!equilibrium}%
.
\end{theorem}

There exists various different proofs of the existence of Nash equilibria for
two-person zero-sum games \cite{Straffin1993, Gibbons1992, Aubin1993,
Neumann1947}. In this note we shall focus on its equivalence with the Dutch
Book Theorem.

The minimax inequality%
\[
\max_{a\in A}\min_{b\in B}g\left(  a,b\right)  \leq\min_{b\in B}\max_{a\in
A}g\left(  a,b\right)
\]
is proved in exactly the same way for surreal payoff functions as for real
payoff functions. The game is said to be in \emph{equilibrium} when these
quantities are equal. The common value is the \emph{value of the game}. For
any mixed strategy $P$ for Alice the minimum of $g\left(  P,Q\right)  $ over
distributions $Q$ is attained when $Q$ is concentrated in a point, i.e.
$Q=\delta_{b}$ for some pure strategy $b\in B.$ Thus%
\begin{equation}
\min_{Q}g\left(  P,Q\right)  =\min_{b}\sum_{a}g\left(  a,b\right)  \cdot
p_{a}. \label{min}%
\end{equation}
To maximize this over all surreal-valued distributions $P$ is a linear
programming problem and can be solved by the same methods as if the payoff
functions were real valued. In particular there exists a surreal valued
distribution that maximizes (\ref{min}). Using this argument we see that
minimax and maximin are obtained even for mixed strategies.

\begin{proof}
[Proof of equivalence of Thm. \ref{DutchBookThm1} and Thm. \ref{ZeroSumThm}%
]Assume that for a two person zero-sum game there exists a value $\lambda$
with optimal strategies $P$ and $Q.$ Then $\lambda<0$ leads to the existence
of a Dutch book and $\lambda\geq0$ leads to the existence of a distribution
$P$ satisfying (\ref{Positiv}).

Assume that the Dutch Book Theorem holds. Assume that there exist a surreal
number $\lambda$ such that
\[
\max_{P}\min_{Q}g\left(  P,Q\right)  <\lambda<\min_{Q}\max_{P}g\left(
P,Q\right)
\]
Consider the payoff function $f\left(  a,b\right)  =g\left(  a,b\right)
-\lambda$. According to the Dutch Book Theorem there exists a probability
distribution $P$ on $A$%
\[
\sum_{a\in A}p_{a}\cdot f\left(  a,b\right)  \geq0
\]
for all $b\in B;$ or there exists a probability distribution $Q$ on $B$ such
that
\[
\sum_{b\in B}q_{b}\cdot f\left(  a,b\right)  <0
\]
for all $a\in A$. Therefore there exists a probability distribution $P$ on $A$
such that%
\begin{equation}
\sum_{a\in A}p_{a}\cdot g\left(  a,b\right)  \leq\lambda\label{mindreend}%
\end{equation}
for all $a\in A$ or there exists a probability distribution $Q$ on $B$ such
that%
\begin{equation}
\sum_{b\in B}q_{b}\cdot g\left(  a,b\right)  \geq\lambda\label{stoerreend}%
\end{equation}
for all strategies $a\in A.$ Inequality (\ref{mindreend}) contradicts that
$\lambda<\min_{Q}\max_{P}g\left(  P,Q\right)  $ and Inequality
(\ref{stoerreend}) contradicts that $\max_{P}\min_{Q}g\left(  P,Q\right)
<\lambda.$ Hence, $\max_{P}\min_{Q}g\left(  P,Q\right)  =\min_{Q}\max
_{P}g\left(  P,Q\right)  .$
\end{proof}

The importance of the proof that the Dutch Book Theorem is equivalent to the
existence of a Nash equilibrium for two-person zero-sum games is that it means
that the two results refer to the same type of rationality. The next example
show that the use of using surreal probabilities may make the difference
between winning and losing.

\begin{example}
\label{ExSur}%
\begin{table}[tbp] \centering
$%
\begin{tabular}
[c]{|l|l|l|}\hline
$g$ & $a_{1}$ & $a_{2}$\\\hline
$b_{1}$ & $1+1/\omega$ & $-1-2/\omega$\\\hline
$b_{2}$ & $-1$ & $1+1/\omega$\\\hline
\end{tabular}
\ \ $\ \caption{Payoff for Alice.\label{Table1}}%
\end{table}%
Consider the payoff function in Table \ref{Table1}. If Alice ignores
infinitesimals her optimal strategy is the distribution $\left(
1/2,1/2\right)  ,$ which gives a payoff function for Bob that is $-1/2\omega$
if $b=b_{1}$ and $1/2\omega$ if $b=b_{2}.$ In this case Bob could win the game
by choosing $b=b_{1}.$ The minimax optimal strategy for Alice is the mixed
strategy $\left(  1/2+\frac{1}{4\left(  \omega+1\right)  },1/2-\frac
{1}{4\left(  \omega+1\right)  }\right)  .$ If she choose this mixed strategy
the payoff is always positive and she will win the game.

One should note that playing this game is not very different from playing the
game where we have scaled the payoff up by a factor $\omega$ (see Table
\ref{Table2}). We may also scale up Bob's optimal strategy by a factor
$4\left(  \omega+1\right)  $ to obtain $\left(  2\omega+3,2\omega+1\right)  .$
Therefore an optimal strategy for Alice is to play the game $4\left(
\omega+1\right)  $ "times" in parallel in such a way that $a_{1}$ is "chosen
$2\omega+3$ times" and $a_{2}$ is "chosen $2\omega+1$ times ".
\begin{table}[tbp] \centering
$%
\begin{tabular}
[c]{|l|l|l|}\hline
$g$ & $a_{1}$ & $a_{2}$\\\hline
$b_{1}$ & $\omega+1$ & $-\omega-2$\\\hline
$b_{2}$ & $-\omega$ & $\omega+1$\\\hline
\end{tabular}
\ \ \ $\ \ \caption{Payoff for Alice multiplied by $\omega .$\label{Table2}}%
\end{table}%

\end{example}

If a two-persons zero-sum game has a Nash equilibrium pair $\left(  \tilde
{a},\tilde{b}\right)  ,$ which is always the case if $A$ and $B$ are finite,
then $\sup_{a\in A}g\left(  a,\tilde{b}\right)  =g\left(  \tilde{a},\tilde
{b}\right)  $ and therefore $\inf_{b\in B}\sup_{a\in A}g\left(  a,b\right)
\leq g\left(  \tilde{a},\tilde{b}\right)  .$ Similarly, $\sup_{a\in A}%
\inf_{b\in B}g\left(  a,b\right)  \geq g\left(  \tilde{a},\tilde{b}\right)  .$
Thus, the game is in equilibrium%
\index{equilibrium}
and the value of the game%
\index{value of game}
is $g\left(  \tilde{a},\tilde{b}\right)  .$ In particular all Nash equilibria
have the same value. The same argument holds for mixed strategies.

\section{Dutch books for short games}

Surreal numbers are totally ordered and never confused with each other. Games
that are not surreal number are confused with a small or large interval of
surreal numbers. For instance $\mathord\ast$ is confused with $0$ and the game
$\left\{  100\mid-100\right\}  $ is confused with any number between $-100$
and $100.$ Before formulating a Dutch Book Theorem for general combinatorial
games we need to introduce the \emph{mean value} $\mu\left(  G\right)  $ of a
short game $G.$ A game $G$ is said to be \emph{short} if it only has finitely
many positions$.$ Our recursive definition of games allows transfinite
recursion and games that are not short, but for the definition of mean values
we shall focus on the short games. Note that if a short game is a number then
it is a dyadic fraction.

The mean value of a game $G$ is a real number $\mu\left(  G\right)  $ that
satisfies the following mean value theorem.

\begin{theorem}
[\cite{Conway1976}]\label{MeanValueThm}If $G$ is a short game then there
exists a natural number $m$ and a number $\mu\left(  G\right)  $ that
satisfies%
\[
n\cdot\mu\left(  G\right)  -m\leq n\cdot G\leq n\cdot\mu\left(  G\right)  +m
\]
for all natural numbers $n.$
\end{theorem}

Mean values of short games can be calculated by the \emph{thermographic
method} described in \cite{Conway1976} and using this method it is easy to see
that the mean value of a short game is always a rational number. Mean values
of games share some important properties with mean values of random variables.
For instance we have

\begin{itemize}
\item $\mu\left(  n\cdot G\right)  =n\cdot\mu\left(  G\right)  ,$

\item $\mu\left(  G+H\right)  =\mu\left(  G\right)  +\mu\left(  H\right)  ,$

\item $G\geq0\Rightarrow\mu\left(  G\right)  \geq0,$

\item $\mu\left(  1\right)  =1.$
\end{itemize}

\begin{example}
\label{ExZeroMean1}The game $G=\left\{  1\mid\left\{  0\mid-2\right\}
\right\}  $ that is illustrated in Figure \ref{mean0}, satisfies $G>0.$ In the
game $n\cdot G$ Right can only play in a sub-game where Left has not played
and the response optimal for Left is always to answer a move of Right by a
move in the same sub-game. From this one sees that $n\cdot G\leq1$ and
therefore that $\mu\left(  G\right)  =0.$ We see that Left may win a game
\emph{for sure} although the game has mean value zero!
\end{example}

The setup is as before that each bookmaker $b\in B$ tells Alice which game he
wants to play if a certain horse $a\in A$ wins. Alice is going to play Left
and the bookmaker or the bookmakers are going to play Right. After certain
bookmakers have been accepted the bookmakers choose natural numbers
$n_{b},b\in B$ and combine these into a super game $\sum_{b\in B}n_{b}\cdot
G\left(  a,b\right)  $ that will depend on which horse wins. We say that we
have a \emph{Dutch book} if there exists natural numbers $n_{1},n_{2}%
,\cdots,n_{k}$ such that Alice will lose the game
\begin{equation}
\sum_{b\in B}n_{b}\cdot G\left(  a,b\right)  \label{dutchbook}%
\end{equation}
for any value of $a.$ Otherwise the set of game valued payoff functions is
said to be \emph{coherent}. If all the games are short surreal numbers then
this notion of coherence is equivalent to the definition of coherence given in
Section \ref{SecSurreal}.

Alice is allowed to choose that the game should be played a number of times in
parallel. With this setup we get the following version of the Dutch Book Theorem.

\begin{theorem}
\label{DutchThm}If a payoff function $G\left(  a,b\right)  ,a\in A,b\in B$
with short games as values, is coherent then either exists a probability
vector $a\rightarrow p_{a}$ and a natural number $n$ such that $np_{a}%
\in\mathbb{N}$ and the game \qquad\
\begin{equation}
\sum_{a}\left(  np_{a}\right)  \cdot G\left(  a,b\right)  >0,\text{ for all
}b\in B, \label{positive mean}%
\end{equation}
or there exist natural numbers $n_{1},n_{2},\cdots,n_{k},$ a natural number
$n$ and a probability vector $a\rightarrow p_{a}$ such that both games
(\ref{dutchbook}) and (\ref{positive mean}) have mean value $0.$
\end{theorem}

\begin{proof}
We apply the existence of an equilibrium in the two-person zero-sum game with
payoff function $\left(  a,b\right)  \rightarrow\mu\left(  G\left(
a,b\right)  \right)  .$ If the value of the two-person zero-sum game is
negative then the game (\ref{dutchbook}) is negative if the coefficients
$n_{1},n_{2},\cdots,n_{k}$ are large enough. If the value of the two-person
zero-sum game is non-negative there exists a probability vector $a\rightarrow
p_{a}$ such that
\[
\sum_{a}p_{a}\cdot\mu\left(  G\left(  a,b\right)  \right)  \geq0.
\]
The mean value of a short game is a rational number. Therefore the probability
vector $a\rightarrow p_{a}$ can be chosen with rational point probabilities.
Hence, there exists a natural number $m$ such that $m\cdot p_{a}$ is an
integer for all $a\in A.$ Therefore%
\begin{align*}
0  &  \leq m\sum_{a}p_{a}\cdot\mu\left(  G\left(  a,b\right)  \right) \\
&  \leq\sum_{a}mp_{a}\cdot\mu\left(  G\left(  a,b\right)  \right) \\
&  =\mu\left(  \sum_{a}mp_{a}\cdot G\left(  a,b\right)  \right)  .
\end{align*}
If
\[
\mu\left(  \sum_{a}mp_{a}\cdot G\left(  a,b\right)  \right)  >0
\]
then there exists a natural number $k$ such that
\[
k\sum_{a}mp_{a}\cdot G\left(  a,b\right)  >0
\]
and the game defined in (\ref{positive mean}) is winning for Alice who plays
as Left when $n=km.$ Otherwise
\begin{equation}
\mu\left(  \sum_{a}mp_{a}\cdot G\left(  a,b\right)  \right)  =0.
\label{meanzero}%
\end{equation}

\end{proof}

Here we should note that our short-game-valued Dutch Book Theorem stated there
are \emph{three }cases that are \emph{not exclusive:}

\begin{enumerate}
\item Dutch book.

\item Positive mean.

\item Zero mean.
\end{enumerate}

As we saw in Example \ref{ExZeroMean1} a game with mean zero may be positive
or negative. Therefore a decision strategy in which only games with positive
means are acceptable will exclude some games that one will win for sure and a
decision strategy where games with non-negative mean are acceptable will
include some games that are lost for sure. The most reasonable solution to
this problem seems to be to accept or reject according to the mean payoff with
respect to some probability distribution, but leave the cases with mean zero
undecided because a more detailed non-probabilistic analysis is needed for
these cases.

\section{More on infinitesimals}

The Dutch Book Theorem for short games only used rational valued mean values.
One may hope for a better Dutch Book Theorem if we allow also allow a mean
value function with infinitesimal surreal numbers as mean values. For short
games this will not solve the problem.

\begin{definition}
A game $G$ is said to be \emph{strongly infinitesimal} if $-s\leq G\leq s$ for
any surreal number $s>0.$
\end{definition}

\begin{example}
The game $\left\{  0\mid\mathord\ast\right\}  $ is called \emph{up} and
denoted $\mathord\uparrow.$ It is easy to check that $\mathord\uparrow>0.$ The
game $\mathord\uparrow$ is infinitesimal (check how Left can win
$2^{-s}-\mathord\uparrow$). One can prove that any infinitesimal short game is
strongly infinitesimal \cite{Schleicher2006}.
\end{example}

An interesting situation is when all games $G\left(  a,b\right)  $ are
infinitesimal. In this case the Dutch Book Theorem for games as formulated in
Theorem \ref{DutchThm} tells exactly nothing because the mean value of
strongly infinitesimal games would always be $0$ even if surreal mean values
are allowed. But if all games are infinitesimal one could shift to a different
"mean value" concept. For short games one compares the game with $n\cdot1$ and
the game 1 can be considered as a unit in the theory. For infinitesimal short
games one can compare with the infinitesimal game $\mathord\uparrow$ instead.
It is possible to define an \emph{atomic mean value} such that
$\mathord\uparrow$ has mean $1,$ but the proofs are more involved. One can
also prove a version of the Dutch Book Theorem for infinitesimally short games
that involves three cases. The three cases are Dutch book, positive mean, and
some games $G$ that cannot be analyzed in the sense that their atomic mean
value is zero. Although infinitesimal games can be treated with their own mean
value concept this will not solve all problems because games that are not
infinitesimal may sometimes be combined into strongly infinitesimal games. A
simple example consist of the games $1$ and $\mathord\uparrow-1$ whose sum is
the strongly infinitesimal game $\mathord\uparrow.$

\section{Discussion}

In any frequency interpretation of probability theory, probabilities should be
interpreted as limits of frequencies. Obviously surreal probabilities cannot
have such interpretations because a frequency interpretation cannot
distinguish between surreal probabilities that have an infinitesimal
difference. This leads us to the following conclusion: frequency probabilities
are real numbers but uncertainty should in general be modelled by convex sets
of surreal numbers.

In a \emph{subjective Bayesian} approach%
\index{Bayesian!approach}
to probability and statistics one will assign probabilities expressing the
individual feeling of how probable or likely an event is. Many subjective
Bayesians justify this point of view by reference to the Dutch Book Theorem.
We note that unlike some of the modification by Savage et al. neither our
formulation of the Dutch Book Theorem nor its original formulation of de
Finetti has any reference to subjectivity. For short-game valued payoffs even
the one-to-one correspondence between probability and coherent decisions
breaks down. Experiments have demonstrated that most people have a bad
intuition of probabilities and are unable to assign probabilities in a
consistent%
\index{consistent}
manner. It should be even harder to make a consistent distinction between the
probabilities $1/3$ and $1/3+1/\omega$ although the Dutch Book Theorem give
the same type of justification for surreal probabilities as for real probabilities.

We have seen that from a mathematical point of view uncertainties may be
modeled by a convex set of surreal probability vectors, but the reader may
wonder why infinitesimals do normally not appear in probability theory.
Actually there are many real numbers that never appear as probabilities. For
instance all the numbers that \emph{do} appear are \emph{computable}, and
there are only countably many computable numbers. Therefore, it seems that the
use of surreal numbers is an idealization that is not worse than the use of
real numbers as subjective probabilities. At the moment two-person zero-sum
games like the ones described in Example \ref{ExSur} are the only known kind
of calculations that gives surreal valued probabilities as results.

In this paper we used the operations $+$ and $\cdot$ to define Dutch books and
coherence. These operations refer to ways of combining games into new games.
It is an open question what kind of Dutch Book Theorem one would get if other
ways of combining games were considered.

For social games with several players and surreal-valued payoff functions we
have not been able to prove existence of a Nash equilibrium, because one
cannot use the usual fixed-point results that rely heavily on the topology of
the real numbers. We shall not discuss it here as it has less interest for our
understanding of what probabilities are.

\section*{Acknowledgements}

Thanks to Wouter Koolen-Wijkstra, Peter Gr\"{u}nwald, and Mogens Esrom Larsen
for many useful comments and discussions.

This work has was supported by the European Pascal Network of Excellence.

\bibliographystyle{plain}
\bibliography{database1}

\end{document}